\documentclass[12pt]{article}

\usepackage{latexsym}

\newtheorem{theorem}{Theorem}[section]
\newtheorem{prop}[theorem]{Proposition}
\newtheorem{lemma}[theorem]{Lemma}

\begin{document}

\title{The Cylinder Theorem in ${\cal H}^2\times R$ }
\author{ J. Lucas M. Barbosa\thanks{This author was supported by CNPq} and M. P. do Carmo \thanks{When this work was essentially ready for publication, Manfredo died.}}
\date{\today}
\maketitle

\abstract
We consider cylinders in ${\cal H}^2\times R$ (see definitions in the introduction) and prove that a complete and connected surface in ${\cal H}^2\times R$ with the vanishing of the Gauss and extrinsic curvatures is a cylinder.

\section{Introduction}

Except for a note by A.V. Pogorelof in the Doklady in 1956 (see details in the Dover edition of Differential Geometry of Curves and Surfaces p. 415), the first time that the theorem of the cylinder in $R^3$ was in print was in a paper by P. Hartman and L. Nirenberg \cite{HN} which treats a more general situation but mentions explicitly the case of surface, namely, if the Gauss curvature vanishes everywhere, the surface is a cylinder. Direct proofs of this particular case were published almost simultaneously, by Massey \cite{Massey} and Stoker \cite{Stoker}.

In this work our ambient space will be the product ${\cal H}^2\times R$ where ${\cal H}^2$ stands for the hyperbolic plane with curvature $-1$.
In this ambient space, we define a {\em cylinder} as the surface given by $\alpha \times R$ where $\alpha $ is any regular curve in ${\cal H}^2$.
As it is well known, a surface in a tree dimension ambient space, besides the Gaussian Curvature, has, at each point, the extrinsic curvature, which is the product of the principal curvatures. In $R^3$ these two curvatures coincide. In ${\cal H}^2\times R$, they are in general different and give different information about the surface.

 The goal of this paper is to prove the following result.

{ \em Let $M$ be a complete and connected surface in ${\cal H}^2 \times R$. Then $M$ is a cylinder if and only if both
the Gaussian curvature and the extrinsic curvature vanish.}

\section{Preliminaries}

Along this work, we are going to use moving frames having as reference the book \cite{oneil}.

\vspace{0,5cm}

In this section, we are going to compute the Gaussian and the extrinsic curvatures of a cylinder.

\vspace{0,5 cm}

Let ${\cal H}^2$ be the hyperbolic space of dimension $2$ and Gaussian curvature $-1$.
Consider the product space ${\cal H}^2\times R$.

\vspace{0,5cm}

As established in the introduction, a {\bf cylinder} in this space is the surface
$M = \gamma \times R$, where $\gamma $ is a curve in ${\cal H}^2$.

\vspace{0,5cm}

To study the geometry of such cylinders we proceed as follows. First, we observe that
the metric of ${\cal H}^2\times R$ is given by
$$
d\sigma^2 = d\zeta^2 + dt^2,
$$
where $d\zeta$ is the metric of ${\cal H}^2$ and $dt$ is the metric of $R$.
The covariant differential $\bar D$ in ${\cal H}^2 \times R$ decomposes naturally as
$$
{\bar D} = D + d,
$$
where $D$ is the covariant differential in ${\cal H}^2$ and $d$ is the standard differentiation in $R$.
If we represent by $\partial /\partial t$ the unit vector field tangent to the lines $p\times R$, where $p$
is any point of ${\cal H}^2$, then we obtain
$$
{\bar D}(\partial/\partial t) = 0.
$$
Let $\gamma $ be a curve in ${\cal H}^2$ that we assume parameterized by the arc length $s$.

Take the following frame field adapted to the cylinder $M =\gamma \times R$: $e_1 = \partial/\partial t$, $e_2 = \gamma'$, $e_3$ normal to $\gamma$ and tangent to ${\cal H}^2$.
The dual forms are then $\omega_1 = dt$, $\omega_2 = ds$ and $\omega_3 = 0$. The metric of $M$ is then given by $dt^2 + ds^2$, hence $M$ is isometric to Euclidean plane and so its Gaussian curvature is identically zero.

Represent by $\omega_{ij}$ the connection forms in $M$, given by ${\bar D}e_i = \sum \omega_{ij}\wedge e_j$, where $\omega_{ij} + \omega_{ji} = 0$. We then have
$$
0 = {\bar D}(\partial/\partial t) = {\bar D}e_1 = \omega_{12}e_2 + \omega_{13}e_3
$$
Hence, $\omega_{12} = 0$ and $\omega_{13}=0$. It follows that
$$
{\bar D}e_2 =\omega_{21}e_1 + \omega_{23}e_3 = \omega_{23}e_3
$$
We have that
$$
{\bar D}e_3 = \omega_{31}e_1 + \omega_{32}e_2
$$
and so, the forms $\omega_{31}$ and $\omega_{32}$ give information about the principal curvatures $k_1$ and $k_2$.
Since $e_2 = \gamma'$ then ${\bar D}e_2 $ must be tangent to ${\cal H}^2$ and must annihilate the vector $\partial/\partial t$.  Therefore, $\omega_{23}$ must be of the form $-\lambda \omega_2$. The fact that $\omega_{31} =0$ and $\omega_{32} = \lambda \omega_2$ tell us that the extrinsic curvature $k_1k_2$ at any point of $M$ is zero. We thus have proved the following proposition that is the necessary condition of the result we want to prove in this work.

\begin{prop}
Let $M$ be a cylinder in ${\cal H}^2 \times R$. Then $M$ has both extrinsic and intrinsic curvatures
identically zero.
\end{prop}

\section{Asymptotic lines  passing through \\ parabolic points of $M$.}

Consider a complete surface $M$ with extrinsic and intrinsic curvatures identically zero in ${\cal H}^2 \times R$. Represent by $P$ the set of planar points of $M$, that is, the points where principal curvatures are zero. Define $U = M - P$, the set of parabolic points of $M$, that is, the points with one of the principal curvatures zero and the other, not zero. The direction with principal curvature zero is called an asymptotic direction and a curve whose tangent vectors are in asymptotic directions is an asymptotic line. It is clear that by a parabolic point passes a unique asymptotic line. Finally, we observe that $P$ is a closed set of $M$ and, consequently, $U$ is an open set.

\begin{prop}
The unique asymptotic line through a parabolic point of $M$ is an open segment of geodesic of the ambient space.
\end{prop}

Take $p \in U$. Let $\gamma$ be the asymptotic line passing through $p$ parameterized by arc length. In an open neighborhood of $p$ take an orthonormal frame field $e_1$, $e_2$, $e_3$ where $e_3$ is normal to $M$ and $e_1 = \gamma'$. Let $\theta_1$, $\theta_2$ be the dual forms corresponding to $e_1$ and $e_2$. The connection forms will be given by
$$
\begin{array} {rcccccl}
{\bar D}e_1 & = &                 &   & \theta_{12} e_ 2 & + & \theta_{13} e_3 \\
{\bar D}e_2 & = & \theta_{21} e_1 & + &                  &   & \theta_{23} e_3 \\
{\bar D}e_3 & = & \theta_{31} e_1 & + & \theta_{32} e_2  &   &                 \\
\end{array}
$$
By the choice made, the fields $e_1$ and $e_2$ are principal being the principal curvatures $k_1=0$ e $k_2\neq 0$.
Hence, ${\bar D}_{e_1}e_3 =0$ while ${\bar D}_{e_2}e_3 = k_2 e_2$. Thus, $\theta_{31} (e_1) =0$, $\theta_{32}(e_1)=0$ and $\theta_{31}(e_2) = 0$, while $\theta_{32} = k_2 \theta_2$. It follows that
$$
\begin{array} {rcccccl}
{\bar D}e_1 & = &                 &   & \theta_{12} e_ 2 &   &      \\
{\bar D}e_2 & = & \theta_{21} e_1 & + &                  &   & -k_2\theta_2 e_3 \\
{\bar D}e_3 & = &                 &   & k_2\theta_2 e_2  &   &                 \\
\end{array}
$$

We can also obtain information about $\theta_{12}$ by observing that, since $\theta_{31} = 0$, then
$$
0 = d\theta_{31} = \theta_{32} \wedge \theta_{21} = k_2 \theta_2 \wedge \theta_{21}
$$
Therefore, $\theta_{21} = \lambda \theta_2$.

Rewriting the above equations in terms of  $\gamma$ we obtain
\begin{eqnarray} \nonumber
\gamma'' & = & 0 \\ \nonumber
{\bar D}_{\gamma'}e_2 & = & 0 \\ \nonumber
{\bar D}_{\gamma'}e_3 & = & 0
\end{eqnarray}
Thus,  $\gamma$ is a geodesic segment in ${\cal H}^2\times R$ contained in $M$. Therefore, our proposition is proved. We also conclude that the fields $e_2$ and $e_3$ are covariantly constants along $\gamma$.

\vspace{2mm}

We can go further and prove that if one of the asymptotic lines is extended as far as possible, it never goes into $P$. We state that in the following proposition. % We observe that it generalizes the well-known Massey Theorem (\cite{Massey}).

\begin{prop} \label{massey}
 Let $M$ be a complete surface of ${\cal H}^2\times R$ with both, extrinsic and intrinsic curvatures identically zero. Let $\gamma$ be a maximal asymptotic line passing through a parabolic point $p \in U \subset M$. Let $P$ be the set of planar points in $M$. Then $\gamma \cap P = \emptyset$.
\end{prop}

To prove this proposition we make use of the following lemma.

\begin{lemma}\label{lema 2}
Let $M$ be a surface in ${\cal H}^2\times R$ Assume that $M$ has both the extrinsic and intrinsic curvatures identically zero. Let $H$ be the mean curvature of $M$. Along an asymptotic line $\gamma$, mentioned above, parameterized by arc length $s$, we have
$$
\frac{d^2}{ds^2}\left ( \frac{1}{H}\right) = 0
$$
\end{lemma}
To prove this lemma we will use the frame field already mentioned.
As we have seen $\theta_{21} = \lambda \theta_2$ and $\theta_{32} = k_2\theta_2$. I claim that
\begin{eqnarray} \nonumber
\lambda' - \lambda^2  & = & 0 \\ \nonumber
k_2' - \lambda k_2 & = & 0
\end{eqnarray}
The proof of these equations is simple. We start from the equation $\theta_{21} = \lambda \theta_2$. Observe that
$d\lambda = \lambda_1 \theta_1 + \lambda_2 \theta_2$. By one side we have $d\theta_{21} = K\theta_1\wedge\theta_2 = 0$. By the other side
$$d\theta_{21} = d\lambda\wedge\theta_2 + \lambda\theta_{21}\wedge \theta_1 = \lambda_1\theta_1\wedge \theta_2 - \lambda^2\theta_1\wedge\theta_2 = (\lambda' - \lambda^2) \theta_1\wedge\theta_2.$$
Therefore, $\lambda' - \lambda^2 = 0$.

The proof of the second equation is done in a similar way. We have $\theta_{32} = k_2\theta_2$. Write down
$dk_2 = (k_2)_1\theta_1 + (k_2)_2 \theta_2$.  By one side we have $d\theta_{32} = \theta_{31}\wedge \theta_{12} = 0$ since $\theta_{31}=0$. By the other side
$$
d\theta_{32} = d (k_2\theta_2) = dk_2\wedge \theta_2 + k_2 \theta_{21}\wedge\theta_1 = (k_2)_1\theta_1\wedge\theta_2 - \lambda k_2 \theta_1\wedge \theta_2
$$
Thus, $(k_2)_1 - \lambda k_2 = 0$ and along $\gamma$ we have $k_2' - \lambda k_2 = 0$. This concludes the deduction of the two equations.

The proof of the lemma can be done in the following way. First we observe that, along  $\gamma$ we have $2H = k_2$. Hence, it is sufficient to prove the result for $k_2$.
$$
\frac{d}{ds} \left( \frac{1}{k_2} \right) = - \frac{k_2'}{k_2^2} = - \frac{\lambda k_2}{k_2^2} = - \frac{\lambda}{k_2}
$$
$$
\frac{d^2}{ds^2}\left( \frac{1}{k_2} \right) = - \frac{d}{ds}\left( \frac{\lambda}{k_2} \right) = - \frac{k_2\lambda' - \lambda(k_2)'}{k_2^2} = -\frac{\lambda' - \lambda^2}{k_2} = 0
$$
This proves the lemma.

{\bf Proof (of the proposition):} Suppose that a maximal asymptotic line  $\gamma$ passing through $p \in U$ contains a point $q\in P$. Since $\gamma$ is connected and $U$ is open in $M$ then there exists a point  $p_0 = \gamma(s_0)$ in $\gamma$ such that $p_0 \in P$ and the points $\gamma(s)$ with $s<s_0$ belong to $U$.

On the other hand, we have seen that, along $\gamma$, for $s<s_0$, we have $H(s) = 1/(as +b)$ where $a$ and $b$ are constants. Since the points of $P$ have mean curvature zero, then we should have
$$
0 = H(p_0) = \lim_{s \rightarrow s_0} H(s) = \lim_{s\rightarrow s_0} \frac{1}{as + b}
$$
that is a contradiction and concludes the proof.

\section{Umbilic surfaces with zero extrinsic and intrinsic curvatures}

Let $M$ be a connected surface in ${\cal H}^2\times R$ which has both, extrinsic and intrinsic, zero curvatures. As before, let $P$ be the set of points of $M$ where the principal curvatures are equal to zero and  $U$ be its complement.

\begin{prop} \label{foliation}
Let $P$ be a complete and connected surface in  ${\cal H}^2\times R$ that has Gaussian and extrinsic curvatures identically zero.  Then, in a neighborhood of any of its points $P$ contains a product of an arc of geodesic of ${\cal H}^2$ by a line segment..
\end{prop}
%\begin{prop}
%$P$ is foliated by lines and any neighborhood of its points contains a product of a geodesic of ${\cal H}^2$ by a segment of line.
%\end{prop}

Before we start the proof of this proposition lets examine how are the geodesics of ${\cal H}^2 \times R$. Let $\gamma $ be one such geodesic and let ${\tilde \gamma}$ be the projection of $\gamma $ on ${\cal H}^2 $. Let $f$ be the projection of $\gamma $ on $R$.
\begin{lemma}
$\gamma$ is a geodesic of ${\cal H}^2 \times R$ if and only if its projection $\tilde \gamma$ on ${\cal H}^2$ is a geodesic in ${\cal H}^2$ and $f(s) = as+b$, being $a$ and $b$ real constants.
\end{lemma}
This lemma is a simple consequence of the fact already mentioned that the covariant differentiation $\bar D$ of the ambient space decomposes
as ${\bar D} = D + d$ where $D$ is the covariant differentiation on ${\cal H}^2$ and $d$ is the standard differentiation on $R$.
Then, if ${\bar \alpha}: (0,1) \rightarrow {\cal H}^2\times R$ is a curve in the ambient space  such that ${\bar \alpha} = \alpha + f$ where $\alpha:(0,1) \rightarrow {\cal H}^2$ and $f:(0,1) \rightarrow R$, we have
$$
\frac{{\bar D} {\bar \alpha}}{ds} = \frac{D\alpha}{ds} + \frac{df}{ds}\frac{\partial}{\partial t}
$$
and
$$
\frac{{\bar D}^2 {\bar \alpha}}{ds^2} = \frac{D^2\alpha}{ds^2} + \frac{d^2f}{ds^2}\frac{\partial}{\partial t}
$$
Thus $\frac{{\bar D}^2 {\bar \alpha}}{ds^2} = 0$ if and only if $\frac{D^2\alpha}{ds^2} = 0$ and  $\frac{d^2f}{ds^2} = 0$ which proves the lemma.

\vspace{0,5cm}

We return now to the proof of the proposition.
Let $P$ be the surface for which $k_1 = k_2 =0$ and the Gaussian curvature is also zero everywhere. Choose an orthonormal frame field
$e_1, e_2, e_3$, such that $e_3$ is normal to $P$. Then $e_1, e_2$ are tangent to $P$ and are principal vector fields. Let $\omega_1, \omega_2$, be the dual forms corresponding to $e_1, e_2$ and let $\omega_{ij}$ be the connection forms. We will then have
$$
\omega_{31} = k_1\omega_1 = 0 \quad \quad \quad \omega_{32} = k_2 \omega_2 =0.
$$
It follows that the field $e_3$ is parallel along $P$. In fact we will have:
$$
{\bar D}e_3 = \omega_{31}e_1 + \omega_{32}e_2 =0.
$$
what proves our claim.

\noindent Define
$$
f = \langle e_3, \frac{\partial} {\partial t}\rangle,
$$
then $df = 0$ and hence $f$ is constant on each connected component of $P$.

By hypothesis $P$ is connected. Assume also that $f \equiv 0$ on $P$. In this case $\partial/\partial t$ is tangent to $P$ at each point. Thus, the line (integral curve of $\partial/\partial t$) is contained in $P$, or it has a segment contained in $P$.  And this happens on each point of $P$. It follows that $P$ contains a product of an arc of curve $\gamma$ of ${\cal H}^2$ by a segment of line.

We can then take $e_1 = \gamma'(s)$ and $e_2=\partial/\partial t$. It follows $\omega_1 = ds$ and $\omega_2 = dt$. As a consequence of this fact, $\omega_{12} = 0$. Hence
$$
\gamma''(s) = {\bar D}_{e_1}e_1 = \omega_{12}(e_1)e_2 + \omega_{13}(e_1)e_3 = 0
$$
Therefore, $\gamma $ is a geodesic on ${\cal H¨}^2$.

\vspace{0,5 cm}

Now, consider the case in which $f$ is equal to a constant $c \neq 0$. Take $p \in P$, then $p = ({\tilde p}, t)$ where ${\tilde p} \in {\cal H}^2$ and $t\in R$. Consider the surface  ${\cal H}^2\times \{t\}$. Then $p \in {\cal H}^2\times \{t\}$. Let $\gamma$ be the curve intersection of  $P$ with ${\cal H}^2\times \{t\}$. Then we have $\gamma(s) = ({\tilde \gamma}(s), t)$. By taking different values of $t$ we obtain a family of curves $\gamma_t(s) = ({\tilde \gamma}_t(s), t)$.

Take an orthonormal frame field such that $e_1 = \gamma_t'$. Then we have
$\langle e_1, \partial/ \partial t\rangle = 0$. Thus
$$
0 = d\langle e_1, \partial/ \partial t\rangle= \langle \omega_{12} e_2 + \omega_{13} e_3, \partial/ \partial t\rangle = \langle e _2,\partial/ \partial t \rangle\omega_{12}
$$
Here we have used that $e_1$ is a principal direction.

If $\langle e _2,\partial/ \partial t \rangle = 0$ then $e_3$ will be parallel to $\partial/\partial t$. Then $e_1$ and $e_2$ are tangent to ${\cal H}^2$ and so $P$ (or an open set of $P$) coincide with the surface ${\cal H}^2\times \{t\}$ for a value of $t$. But then, in an open set the value of the Gaussian curvature of $P$ is $-1$ and not zero. This is a contradiction.

Therefore we must have $\omega_{12} = 0$.
 Take the function $g = \langle e _2,\partial/ \partial t \rangle$. Observe that
$$
dg = \langle \omega_{21}e_1 + \omega_{23}e_3, \partial/ \partial t\rangle = 0.
$$
Therefore, $g$ is constant and different from zero. Let $\alpha$ be the integral curve of the field $e_2$. The same argument used in the case of the curve $\gamma$ shows that $\alpha$ is a geodesic and that its projection in ${\cal H}^2$ is a geodesic on ${\cal H}^2$. Since the fields $e_1$ and $e_2$ are orthonormal we have that $\gamma$ and $\alpha$ are perpendicular in their intersection. Since $\gamma = ({\tilde \gamma},t) $ and $\alpha = ({\tilde \alpha}, as + b)$  then $\gamma' = ({\tilde \gamma}',0) $ and $\alpha' = ({\tilde \alpha}', a)$ it follows that the geodesics $\tilde \gamma$ and $\tilde \alpha$ are also perpendicular.

Observe now that the construction of these two geodesics can be done at any point of $P$ thus generating a network of perpendicular geodesics on ${\cal H}^2$. But then we would have on ${\cal H}^2$ families of rectangles, that is, geometric figures formed by four arcs of geodesics, that, when they have intersection, they are perpendicular. This can not occur. The result follows from classical hyperbolic geometry or can be obtained as a simple application of Gauss Bonnet Theorem.

This proves that, in a neighborhood of any of its points $P$ contains a product of an arc of geodesic of ${\cal H}^2$ by a line segment. This concludes the proof of the proposition.

\section{The main result}

In this section we are going to prove the sufficient condition of the result we proposed to prove in this work. We start by proving the following proposition that presents a more precise statement of the previous one.

\begin{prop}
Let $M$ be a surface complete and connected in ${\cal H}^2\times R$ that has both, the Gaussian and the extrinsic, curvatures identically zero. Let $P$ be the set of planar points of $M$ and $U = M-P$. Then, the interior of each connected component of $P$, say $P_0$, will be a product $\gamma\times R$ where $\gamma:(a,b)\rightarrow {\cal H}^2 $ is an open arc of geodesic of ${\cal H}^2$. Furthermore $\partial P_0$ are the two lines $\gamma(a) \times R$ and $\gamma(b)\times R$.
\end{prop}

\noindent {\bf Proof:} Let $P_0$ be a connected component of $P$ and $p$ be a point of interior of $P_0$. By proposition (\ref{foliation}) we know there is an open arc of geodesic $\alpha$ and an open interval $I$ such that $p \in \alpha \times I \subset \mbox{interior of } P_0$. We now extend $I$ to a complete line $L$. Observe that $L$ is a geod\'esic of the ambient space and it is also a geod\'esic of $M$ since $L\cap M \supset $ an interval.

I first claim that $L$ cannot leave $P_0$.

If $L$ leaves $P_0$ it will reach a point $q$ of $U$. In a neighborhood $\cal V$ of $q$ we may take a frame field $e_1 = \partial/\partial t$, $e_2$, $e_3$ normal to $M$. This choice is possible since $L$ has an arc contained in $U$ passing through $q$. For this choice we will have
${\bar D}e_1 = 0$ what implie that $\theta_{13}=0$ and so, the arc of $L$ in $U$ is an asymptotic curve. Such curves are unique passing through each point of $U$ and do not intercept $P$, as we have already seen.  This proves the claim.

Our second claim is that $L$ cannot reach any point of $\partial P_0$.

To show this claim we first observe that $L$ is not contained in $\partial P_0$ since $p\in L$ and $p\in \mbox{interior of }P_0$. So, starting from $p$, there is a first point $q\in \partial P_0$ where the line $L$ touches $\partial P_0$.

Take a small neighborhood $\cal V$ of $p$. The lines that cross this neighborhood are separated in two classes, the ones at the right and the ones at the left of $L$. Even if the classification is done using $\cal V$ it can be extended to all points of the line $L$ and the lines will not change class along this procedure. In particular, at the point $q$ we must have lines at the left side and at the right side of L. But, since $q$ belongs to $\partial P_0$ then we will have only lines at one side of $L$. The other ones will contain points of U. This is a contradiction that proves our claim.

It then follows that $P_0$ is foliated by entire lines and furthermore  interior of $P_0$ is the product $\gamma \times R$, where $\gamma$ is an arc of geodesic of ${\cal H}^2$, say $\gamma:(a,b)\rightarrow {\cal H}^2$. Observe that $\partial P_0$ must be a line, as the limit of lines, and, as such, $P_0 = \gamma \times R$, where $\gamma:[a,b]\rightarrow {\cal H}^2$ is a geodesic with extremities, and $\partial P_0$ consists of two lines $\gamma(a)\times R$ and $\gamma(b)\times R$.

This proves our proposition.

\vspace{0,5 cm}

Now we prove our main theorem. As we already said this is the sufficient part of the result we mention in the Introduction whose proof is the goal of this paper. We observe that the necessary condition was already proved in the section "Preliminaries".

 \begin{theorem} \label{main} Let $M$ be a complete and connected surface in ${\cal H}^2 \times R$ for which both
the Gaussian curvature and the extrinsic curvature vanish. Then $M$ is a cylinder.
\end{theorem}

{\bf Proof:} Let $M \subset {\cal H}^2\times R$ be a complete and connected  surface that has the Gaussian and the extrinsic curvatures identically zero. Let $P$ be the set of planar points of $M$ and $U = M -P$ be the set of parabolic points of $M$. Assume $P \ne \emptyset$. By the previous proposition each connected component $P_i$ of $P$ is a product $\gamma_i \times R$ where $\gamma_i:[a_i, b_i]\rightarrow {\cal H}^2$ is an arc of geodesic of ${\cal H}^2$ with extremities, and $\partial P_i$ is formed by two lines $\gamma_i(a_i)\times R$ and $\gamma_i(b_i)\times R$.

Consider the intersection $P_i \bigcap({\cal H}^2\times \{t_0\})$. This is a curve $\alpha$ that has a point $\alpha(s_0)$ on $\partial P_i$. We may extend such a curve getting into U. Through each point of this curve in $U$ it passes a unique asymptotic curve $\gamma$ that is also a geodesic of the ambient space. Such geodesics are either one line or have the form $\gamma(s) = ({\tilde \gamma}(s), as+b)$ where $\tilde \gamma$ is a nontrivial geodesic of ${\cal H}^2$. Assume this is the case.  The family of such geodesics indexed by the points of $\alpha$ in $U$, has a limit that must be the line ${\tilde p}_i\times R$ contained into $\partial P_i$ or a line ${\tilde p}\times R$ not necessarily in $\partial P$.  But then, the limit of the family of the geodesics $\tilde \gamma$ must be the point ${\tilde p}_i$, what is not possible. Therefore, the curves $\gamma $ mentioned  are lines and each connected component of $U$ must be the product of an arc of curve in ${\cal H}^2$ by $R$. This proves the theorem for the case $P\ne \emptyset$.

%\vspace{5mm}

Assume now that $P = \emptyset$. Then $M = U$. Take $p = ({\tilde p}, a) \in M$. Consider the surface $S = {\cal H}^2 \times \{a\}$. It is clear that $M\cap S$ does not contains any open set, otherwise, in this open set we would have that the Gaussian curvature of $M$ would be $-1$ what is note the case.

Take then the curve $\alpha = S\cap M$. Through each point $\alpha(r)$ there is a unique geodesic $\gamma_r$ of the ambient space. Such geodesic is either a line or $\gamma_r(s) = ({\tilde \gamma }_r(s), a_rs +b_r)$, being ${\tilde \gamma}_r$ a geodesic of ${\cal H}^2$.

%It cannot occur that one of the geodesics is a line and the other one is not. In this case, the set of the geodesics that are not lines would %have a limit that would be a line what cannot occur otherwise the projections on ${\cal H}^2$ would be a family of geodesics converging to one %point, what does not exist. This completes the proof of the claim.

%If, for each point $\alpha(r)$, the asymptotic line passing through it is a line then $M$ contains a product of ${\tilde \alpha}$ times a line.

Claim: If for a point of $\alpha$ passes a line contained in $M$ then each geodesic $\gamma_r$ would be a line and $M$ would be a product ${\tilde \alpha}$ times a line.

 Indeed, If such is not true, then there exist points $\alpha(r)$ for witch the corresponding lines of curvature are of the form $\gamma(s) = ({\tilde \gamma }(s), a(s))$. The family of such geodesics would have to converge to one line. But then the geodesics $\tilde \gamma$ would have to converge to one point. This is a contradiction.

Thus our conclusion is that either $M$ is the product $\alpha \times R$ or the lines of curvature passing through points of $\alpha$ are of the form $\gamma(s) = ({\tilde \gamma}(s), a(s))$. Let's prove that the second possibility cannot occur.

\noindent {\bf Claim:} If $\alpha(r_1) \ne \alpha(r_2)$ then the corresponding mentioned geodesics do not intersect each other.

Indeed. We have $\alpha(r_i) = ({\tilde \alpha}(r_i), a)$, $i=1,2$. If the geodesics are lines, then the lines go through two distinct points of ${\cal H}^2$ and so, the do not intersect. Assume that they intersect. The point of intersection is a point of $U = M$ and so, through it passes a unique asymptotic line that is a geodesic of the ambient space. Since the two geodesics are asymptotic lines they would have to coincide. This is a contradiction that proves the claim.

Now we observe that, since $M$ is complete, connected and has zero Gaussian curvature, then $\exp _p:T_pM \rightarrow M$ is an isometry of the Euclidean space $T_pM$ onto $M$. This is in the literature, but can be easily found in \cite{manfredo} page 464 (Corollary and Proposition 7). Then the inverse image of the geodesics $\gamma$ going through points of $\alpha$ are lines on $T_pM$. Since distinct geodesics do not intercept then their inverse image are parallel lines. We can say more, since parallel lines are equidistant, then the geodesics mentioned are also equidistant.

Remember that the metric on ${\cal H}^2\times R$ is given by
$$
d\sigma^2 = d\varsigma^2 + dt^2
$$
where $d\varsigma$ is the metric on ${\cal H}^2$ and $dt$ is the usual metric on $R$, and the metric $d_M$ of $M$ is just the restriction of $d_{\sigma}$ to $M$. Take any two mentioned geodesics, say $\gamma_1$ and $\gamma_2$. Since such geodesics are equidistant, given any point $p_1 \in \gamma_1$ there exists a point $p_2\in \gamma_2$ such that
$$
d_M(\gamma_1, \gamma_2) = d_M(p_1,p_2).
$$
We know that $p_1= ({\tilde p}_1, a_1)$ and $p_2 = ({\tilde p}_2, a_2)$ where ${\tilde p}_1 \in {\tilde \gamma}_1$ and ${\tilde p}_2 \in {\tilde \gamma}_2$ being ${\tilde \gamma}_1$ and ${\tilde \gamma}_2$ geodesics of ${\cal H}^2$. Then we have
$$
d_M(p_1,p_2) \ge d_{\sigma}(p_1, p_2) \ge d_{\varsigma}({\tilde p}_1, {\tilde p}_2).
$$
%We must have $d_M(\gamma_1, \gamma _2) = c$, where $c$ is a real positive constant.
%But
%$$
%d_{\sigma}(p_1,p_2) \ge d_{\varsigma}({\tilde p}_1, {\tilde p}_2) \ge d_{\varsigma}({\tilde p}_1, {\tilde \gamma}_2)
%$$
But, in ${\cal H}^2$, $d_{\varsigma}({\tilde p}_1,{\tilde p}_2)$ is unbounded, that is, for any given positive number $\lambda$ there are points ${\tilde p}_1 \in {\tilde \gamma}_1$ and ${\tilde p}_2 \in {\tilde \gamma_2}$ such that $d_{\varsigma}({\tilde p}_1, {\tilde p}_2)= \lambda$.  This is a contradiction with the fact that $\gamma_1$ and $\gamma_2$ are equidistant and proves the theorem.

%%\section{A note about Manfredo}
%When this work was essentially ready for publication, Manfredo died. The authors discussed issues of this work about 12 days before his death. %The work was completed by the first author and submitted for publication.

%\vspace{5 mm}

\vspace{2cm}

\noindent Jo\~ao Lucas Marques Barbosa \\
Rua Carolina Sucupira 723 Ap 2002 \\
60140-120 Fortaleza - Ce \\
Brazil \\
joaolucasbarbosa@terra.com.br

\vspace{5mm}

\noindent Manfredo Perdig\~ao do Carmo \\
Instituto Nacional de Matem\'atica Pura e Aplicada - IMPA \\
Estrada Dona Castorina 110 \\
22460-320 Rio de Janeiro - RJ \\
Brazil \\
manfredo@impa.br

\end{document}